\documentclass[reqno,  12pt]{amsart}
\usepackage{amsfonts}
\usepackage{amssymb}
\usepackage{amsthm}
\usepackage{upref}
\usepackage{enumerate}

\usepackage{amscd}
\usepackage{bbm}

\usepackage [latin1]{inputenc}

\usepackage{amsmath}
\usepackage{enumerate}

\makeatletter
\def\LaTeX{\leavevmode L\raise.42ex
    \hbox{\kern-.3em\size{\sf@size}{0pt}\selectfont A}\kern-.15em\TeX}
 \makeatother
 
\numberwithin{equation}{section}

\newtheorem{lemma}{Lemma}[section]
\newtheorem{theorem}[lemma]{Theorem} 
\newtheorem{corollary}[lemma]{Corollary}
\newtheorem{proposition}[lemma]{Proposition}

\theoremstyle{definition}

\newtheorem{example}[lemma]{Example}

\newtheorem{assumption}[lemma]{Assumption}

\newtheorem{remark}[lemma]{Remark}

  \newcommand{\e}{\eqref}

\newcommand{\q}{\quad}

\newcommand{\ov}{\overline}
\newcommand{\wt}{\widetilde}

\newcommand{\ti}{\tilde}
\newcommand{\wh}{\widehat}

\renewcommand{\d}{\delta}

   \newcommand{\sgn}{\operatorname{sgn}}

\newenvironment{pf}{\begin{proof}}{\end{proof}}

\def\qqq{\mathrel{\subset\mkern-15mu\lower.38ex\hbox{${\scriptscriptstyle\rightarrow}$}}}

\let\cal\mathcal

\let\Bbb\mathbb

\begin{document}

 
 \title[Universal relations in asymptotic formulas]{Universal relations in asymptotic formulas for orthogonal polynomials}
\author{D. R. Yafaev}
\address{Univ  Rennes, CNRS, IRMAR-UMR 6625, F-35000
    Rennes, France and SPGU, Univ. Nab. 7/9, Saint Petersburg, 199034 Russia}
\email{yafaev@univ-rennes1.fr}
 

\thanks {Supported by RFBR grant No.   17-01-00668 A}


\begin{abstract}
Orthogonal polynomials $P_{n}(\lambda)$ are oscillating functions of $n$ as $n\to\infty$ for $\lambda$ in the absolutely continuous spectrum of the corresponding Jacobi operator $J$.  We show that, irrespective of any specific assumptions on coefficients of the operator $J$, amplitude and phase factors in asymptotic formulas for $P_{n}(\lambda)$  are linked by certain universal relations found in the paper.

Our approach relies on a study of operators diagonalizing
Jacobi operators.  Diagonalizing operators are constructed in terms of orthogonal polynomials $P_{n}(\lambda)$.  They
act from the space $L^2 (\Bbb R)$ of functions into 
the space $\ell^2 ({\Bbb Z}_{+})$ of sequences.   We consider such operators in a rather general  setting and find necessary and sufficient conditions of their boundedness.
   \end{abstract}

   \dedicatory{To  the memory of
    Michail Zakharovitch Solomyak on the occasion of his 90th anniversary}

      \maketitle
  
 


\section{Introduction}

\subsection{Main result}


Jacobi operators $J$ acting in the space $\ell^2({\Bbb Z}_{+})$ are diagonalized (see, e.g.,  the book \cite{AKH}) by unitary operators ${\bf U}: L^2({\Bbb R}; d\rho)  \to  \ell^2({\Bbb Z}_{+})$  ($d\rho$ is the spectral measure of $J$) constructed in terms of the corresponding orthonormal polynomials  $P_{n} (\lambda)$:
 \begin{equation}
({\bf U}  {\bf f} )_{n}=   \int_{\Bbb R} P_{n} (\lambda ){\bf f}(\lambda) d\rho(\lambda),\q n\in {\Bbb Z}_{+}=\{0,1,\ldots\} , \q {\bf f}\in L^2({\Bbb R}; d\rho).
\label{eq:UD}\end{equation}
We are interested in an asymptotic behavior of the  polynomials  $P_{n} (\lambda)$ as $n\to\infty$ for   $\lambda$ in    the absolutely continuous spectrum of $J$. Typically, $P_{n} (\lambda)$ are oscillating functions such that
 \begin{equation}
P_{n}  (\lambda) = 2 \kappa   (\lambda)n^{-r} \cos \big(\omega (\lambda) n^s +\Phi_{n} (\lambda)\big)+ 
O(n^{-\d})  
\label{eq:UD1}\end{equation}
where $s\in (0, 1]$, $r\geq 0$, $\d>1/2$, $\omega'(\lambda) > 0$ and $\Phi_{n}' (\lambda) = o(n^s)$ as $n\to\infty$.
   A more general asymptotic formula of this type is given in Assumption~\ref{ASY}. 

Our main goal is to show  that, irrespective of specific assumptions about matrix elements of the operator $J$, or about the corresponding spectral measure   
\begin{equation} 
d\rho(\lambda)= \tau(\lambda) d\lambda, \q  \tau(\lambda) >0, 
\label{eq:Lag2}\end{equation}
the asymptotic coefficients in formula \e{eq:UD1} are linked by the universal relations 
 \begin{equation} 
  \boxed{
2 r+ s=1}
\label{eq:HXG}\end{equation}
and
\begin{equation} 
  \boxed{
2 \pi \tau(\lambda) \kappa^2 (\lambda) = s   \omega' (\lambda) .}
\label{eq:HXG1}\end{equation}
Given unitarity of the operators \e{eq:UD},
we deduce these relations from asymptotic formulas for time dependent evolution $e^{-i\Theta (J) t} f$  as $t\to\infty$ for suitable functions $\Theta(\lambda)$.  These results were announced in \cite{Y/Lag}; their proofs are given in Sect.~4.


\subsection{Boundedness of semidiscrete Fourier transforms}

Unitarity (and  hence boundedness)  of  operators \e{eq:UD}  are consequences of the 
orthogonal polynomials theory.  In Sect.~3, we   discuss sufficiently arbitrary semidiscrete (that is, acting from a space of functions into a space of sequences) Fourier  operators $V$  imitating their structure
and find   conditions of their  boundedness by direct tools. A typical result about such operators is stated as Theorem~\ref{HL}.   Below $C$ and $c$ are various positive constants whose precise values are irrelevant.

\begin{theorem}\label{HL}
Let an operator  
      $V:  L^2({\Bbb R})\to  \ell^2({\Bbb Z}_{+})$ be defined by the formula  
\begin{equation}
(Vf)_{n}= v_{n} \int_{-\infty}^\infty   e^{i x_{n} \lambda }w (\lambda)f(\lambda) d\lambda
\label{eq:FF2}\end{equation}
where $v_{n}\in {\Bbb R}$, $x_{n}  \in {\Bbb R}$, $x_{n}< x_{n+1}$ and $x_{n}\to\infty$.  If
\begin{equation}
\sup_{R >0}\sum_{n : x_{n}\in (R,R+1)} v_{n}^2 <\infty
\label{eq:discr}\end{equation}
and 
\begin{equation}
|w (\lambda)| \leq C (1+|\lambda|)^{-1},
\label{eq:FF3}\end{equation}
then the operator  $ V:  L^2 ({\Bbb R} ) \to \ell^2 ({\Bbb Z}_{+})$    is bounded. Conversely, if $V$  is bounded for $w$ being the characteristic function of some interval, then condition \e{eq:discr} is satisfied. 
 \end{theorem}
 
 Note that  \e{eq:discr} holds true, if
 \begin{equation}
 x_{n}=    n^{s}, \q s\in (0,1] \q \mbox{and} \q  |v_{n}| \leq C (1+n)^{- r} \q \mbox{with}\q 2r\geq 1-s .
\label{eq:rs}\end{equation}
 
 Let us mention two particular cases of Theorem~\ref{HL}.  Suppose that $w$ is a bounded function with compact support.  If $s=1$, $r=0$, then the boundedness of the operator $V$ 
follows  from the Parseval identity for the exponentials $e^{in\lambda}$. The case $s=1/2$, $r=1/4$ 
can be deduced (see Sect.~2.4)  from   classical results on Hermite (or Laguerre) polynomials. 
It is however clear that 
 
 $1^0$ such  statements should be true for sufficiently arbitrary sequences $x_{n}$ and under less stringent assumptions on the function $w$
 
$2^0$  the theory of orthogonal polynomials is irrelevant for a study of operators \e{eq:FF2}.


Our goal in Sect.~3 is to justify these conjectures.  It is convenient to state the problem  in a more general setting.
 Let $\Phi$,
 \begin{equation}
(\Phi f)(x)= (2\pi)^{-d/2}
\int_{{\Bbb R}^d} e^{ix\xi }f(\xi)d\xi=: \hat{f} (x)
\label{eq:Four}\end{equation} 
be the (adjoint) Fourier transform, and let $W$ be the operator of multiplication by a bounded function $w(\xi)$.  We find in Theorem~\ref{Four} generalizing  Theorem~\ref{HL}  necessary and sufficient conditions of  boundedness of the operator
 \begin{equation}
A= \Phi W: L^2 ({\Bbb R}^d )\to L^2 ({\Bbb R}^d; d M)
\label{eq:FourF}\end{equation} 
where $d M (x)$ is some measure on ${\Bbb R}^d$. This is of course a well studied problem; nevertheless, Theorem~\ref{Four}  seems to be new. Note that the boundedness of operator $A$ for $w(\xi)= (1+ |\xi|)^{-\ell}$ is equivalent to the embedding of the Sobolev space ${\sf H}^\ell ({\Bbb R}^d)$  into $L^2 ({\Bbb R}^d; d M)$.

Let us however mention an important relevant result.  Obviously,  the operator $v\Phi w$ is bounded in the space $L^2({\Bbb R}^d)$  if $v,w \in L^\infty ({\Bbb R}^d)$.  The case of functions $v$ and $w$ with local singularities was studied in the paper \cite{BKS}  by M. Sh. Birman, G. E. Karadzhov and M. Z. Solomyak. The conditions of boundedness of such operators $v\Phi w$ generalize classical Hardy inequalities.


\section{Jacobi operators and orthogonal polynomials}
 
    \subsection{Basic  facts }
      Let us briefly recall necessary results about
    Jacobi operators $J$. They act in the space $\ell^2 ({\Bbb Z}_{+})$ by the formula 
       \[
(J u)_{n}= a_{n-1} u_{n-1}  + b_{n} u_{n }   + a_{n} u_{n+1} ,  \q n\in{\Bbb Z}_{+}, \q a_{-1}=0, 
\]
where the sequences $a_{n }>0$ and $b_{n}=\bar{b}_{n}$ are given  and $u =\{ u_{n}  \}_{n=0}^\infty$.  A minimal symmetric  operator $J_{0}$  is  
defined on a set  of vectors   with only a  finite number of non-zero components $u_{n} $.  The spectra of all self-adjoint extensions $J$ of $J_{0}$  are simple with $e_{0} = (1,0,0,\ldots)^{\top}$ being a generating vector. The   spectral measure of $J$  is defined by the relation $d\rho_{J} (\lambda)=d(E_{J} (\lambda)e_{0}, e_{0})$ where  $E_{J} (\lambda)$      is the spectral family of the operator $J$.
If the   Carleman condition 
 \begin{equation}
    \sum_{n=0}^\infty a_{n}^{-1}=\infty
\label{eq:Carl}\end{equation} 
is satisfied, then the operator $J_{0}$ is essentially self-adjoint.

Orthonormal polynomials  $P_{n}( z )$  can be defined as solutions of the difference equation
 \begin{equation}
 a_{n-1} u_{n-1} (z) +b_{n} u_{n } (z) + a_{n} u_{n+1} (z)= z u_n (z), \q n\in{\Bbb Z}_{+}, 
 \q z\in{\Bbb C}, 
\label{eq:Py}\end{equation}
satisfying the boundary condition  $P_0 (z) =1$. Clearly,  $P_{n } (z)$ is a polynomial of degree $n$:
$P_{n} (z)= p_{n}z^n +\cdots$ with $p_{n} > 0$. It can be shown that
 these  polynomials    are orthogonal and normalized  in the spaces $L^2 ({\Bbb R};d\rho )$: 
    \[
\int_{-\infty}^\infty P_{n}(\lambda) P_{m}(\lambda) d\rho (\lambda) =\d_{n,m};
\]
as usual, $\d_{n,n}=1$ and $\d_{n,m}=0$ for $n\neq m$.  Here $d\rho=d\rho_{J}$ is the spectral measure of an arbitrary self-adjoint extension $J$ of the minimal operator $J_{0}$.
Alternatively,   given a probability measure $d\rho(\lambda)$, the polynomials $P_0 (\lambda),P_1 (\lambda),\ldots,  P_{n} (\lambda),\ldots$ can be obtained by the Gram-Schmidt orthonormalization of the monomials $1,\lambda,\ldots,\lambda^n,\ldots$ in the space $L^2({\Bbb R}_{+}; d\rho)$. 


Recall that   mapping \e{eq:UD}  is unitary, that is,
   \[
{\bf U}^*  {\bf U} =I,\q  {\bf U} {\bf U}^*   =  I
 \]
 (we denote by $I$  the identity operator in different spaces)
  and enjoys the intertwining  property 
$({\bf U}^* J u)(\lambda)= \lambda ({\bf U}^*  u)(\lambda)$.

Let $\Lambda$ be an interval of the absolutely continuous spectrum of the operator $J$ where the spectral measure is given by relation        \e{eq:Lag2}.
On the subspace $L^2 (\Lambda; d\rho)$, mapping \e{eq:UD} can be reduced (if one sets $  \sqrt{\tau(\lambda)} {\bf f}(\lambda)=  f(\lambda)$) to the operator $U=U_{\Delta}: L^2 (\Lambda)\to \ell^2 ({\Bbb Z}_{+})$ defined by the relation
    \begin{equation}
(U f)_{n}=   \int_\Lambda P_{n} (\lambda )\sqrt{\tau(\lambda)} f(\lambda) d\lambda,\q n\in {\Bbb Z}_{+},\q f\in L^2 (\Lambda).
\label{eq:UFf}\end{equation}
Note that the operator $U$ is isometric on the space $L^2 (\Lambda)$.

Below we discuss asymptotic formula  \e{eq:UD1} for $\lambda\in \Lambda$. All our estimates are uniform in $\lambda$ from compact subintervals of an open interval $\Lambda$.

   \subsection{Classical polynomials} 

The Jacobi polynomials $P_{n} (\lambda )= P_{n}^{\alpha,\beta} (\lambda )$ are parametrized by two indices $\alpha>-1$ and $\beta>-1$. The corresponding 
    Jacobi operator  $J =J_{\alpha,\beta}$ has   the  absolutely continuous spectrum $[-1,1]$ with the spectral measure $ d\rho(\lambda) $  defined by \e{eq:Lag2} 
where  the weight function
 \begin{equation}
 \tau (\lambda)=   k (1-\lambda)^{\alpha} (1+\lambda)^{\beta},  \q \alpha ,\beta>-1\q   \lambda\in \Lambda= (-1,1).
\label{eq:Jac1}\end{equation}
 The constant 
$  k $
 is chosen in such a way that measure \e{eq:Lag2} is normalized, i.e., 
$\rho ({\Bbb R})=\rho ((-1,1))=1$. 
Explicit expressions for
 matrix elements  $a_{n}, b_{n}$ of  the Jacobi operator  $J $ can be found, for example, in the books \cite{BE, Sz}, but we do not need them. We only note that
$ a_{n}  =1/2 +  O( n^{-2})$,
$b_{n}  = O( n^{-2} )$.
According to formula (8.21.10)   in the book~\cite{Sz} {\bf the Jacobi polynomials} satisfy relation \e{eq:UD1} where $s=1$,  $r=0$,  $\d=1$ and $\omega(\lambda)=\arcsin\lambda$, 
\[
\kappa(\lambda)= (2 \pi k)^{-1/2}(1-\lambda)^{-(1+2\alpha)/4}(1+\lambda)^{-(1+2\beta)/4} ,
\]
\[
\Phi_{n}(\lambda)= 2^{-1} (\alpha+\beta+1) \arcsin\lambda - \pi(2n+\beta-\alpha) /4.
\]

 {\bf The Laguerre polynomials }  $P_{n} (\lambda )= P_{n}^{\alpha} (\lambda )$  where the parameter $\alpha>-1$ are determined by  recurrence coefficients 
 \begin{equation}
    a_{n} =   \sqrt{(n+1)(n+1+\alpha)} \q\mbox{and}\q     b_{n} =    2n+\alpha+1. 
\label{eq:Lag}\end{equation}
The corresponding Jacobi operators   $J$      have absolutely continuous  spectra
coinciding with $[0,\infty)$ and the weight functions equal
$
\tau (\lambda)=\Gamma(\alpha+1)^{-1}\lambda^\alpha e^{-\lambda}$ where $ \lambda\in {\Bbb R}_{+}$.
According to formula (10.15.1) in \cite{BE}  the Laguerre polynomials satisfy  all  relation \e{eq:UD1} where $s=1/2$,  $r=1/4$,  $\d=3/4$ and $\omega(\lambda)=2\sqrt{\lambda}$, 
\[
\kappa(\lambda)= 2^{-1}\pi^{-1/2} \sqrt{\Gamma(1+\alpha)} \,  \lambda^{-\alpha/2-1/4}  e^{\lambda/2} ,
\q 
\Phi_{n}(\lambda)= \pi n- \pi(2 \alpha+1) /4.
\]

{\bf The Hermite polynomials}   are determined by the Jacobi coefficients 
 $      a_{n}=\sqrt{ (n+1)/2} , \q b_{n}=0$. 
The corresponding Jacobi operator   $J$      has  the absolutely continuous  spectrum
coinciding with the whole axis $\Bbb R$ and the weight function equals
$ \tau (\lambda)=   \pi^{-1/2} e^{-\lambda^2}, \q \lambda\in {\Bbb R}$.
According to formula (10.15.18) in \cite{BE} the Hermite polynomials $P_{n}  (\lambda)$ satisfy, for all  $\lambda\in \Bbb R$, an asymptotic  relation
\begin{equation}
P_{n}  (\lambda)=
   2^{1/2}
  \pi^{-1/4} e^{\lambda^2/2}  (2n+1)^{-1/4} \cos \big( \sqrt{2n+1} \, \lambda-  \pi  n/2\big) +O(n^{-3/4})
\label{eq:H2}\end{equation}
as $n\to\infty$.
We now have $s=1/2$,  $r=1/4$,   $\omega(\lambda)=\sqrt{2\lambda}$ 
and  $\kappa(\lambda)= 2^{-3/4} \pi^{-1/4} e^{\lambda^2/2}$.

  Thus equality  \e{eq:HXG}  is obviously true for all classical polynomials. A direct calculation shows that relation \e{eq:HXG1} is also satisfied.
  
    \begin{remark}\label{Hermite}
  Theorem~\ref{HL} for the particular case $s=1/2$, $r=1/4$ and compactly supported functions $w (\lambda)$   can be deduced from boundedness of  operator \e{eq:UFf}  for the Hermite polynomials.
  Observe first that the operator with kernel $\epsilon_{n}  (\lambda) w(\lambda)  $ is Hilbert-Schmidt if $\epsilon_{n}  (\lambda)= O(n^{-\d} ) $ with $\d> 1/2$. Therefore it follows from formula  \e{eq:H2} that the operator with kernel $  (n+1/2)^{-1/4} \cos \big( \sqrt{n+1/2} \, \lambda-  \pi  n/2\big) w(\lambda)  $  is bounded.  Choosing $n=2m$  or $n=2m+1$ where $m\in{\Bbb Z}_{+}$ and neglecting Hilbert-Schmidt operators, we see that the operators with kernels $  (m+1)^{-1/4} \cos \big( \sqrt{m+1} \, \lambda \big) w(\lambda)  $ and  $  (m+1)^{-1/4} \sin \big( \sqrt{m+1} \, \lambda \big) w(\lambda)  $
are bounded. This ensures boundedness of the operator with kernel $  (m+1)^{-1/4} \exp \big( i \sqrt{m+1} \, \lambda \big) w(\lambda)  $.
  \end{remark}

      \subsection{General polynomials }

  Let us now mention some results for general orthogonal polynomials.  In all the papers discussed below, formula
\e{eq:UD1}  holds true although estimates of the remainders are not always as good as required in \e{eq:UD1}.

  First, we recall the classical result of S.~Bernstein \cite{Bern} (see also     formula (12.1.8) in the G.~Szeg\"o book
  \cite{Sz}). It  states that under some (rather stringent) assumptions on the absolutely continuous spectral measure $d\rho(\lambda)$  supported by $[ -1,1]$, 
  formula \e{eq:UD1} is true with $r=0$, $s=1$, $\omega(\lambda)=\arccos \lambda$ and the amplitude factor $\kappa(\lambda)$ defined by  \e{eq:HXG1}. The phase shift $\Phi_{n}(\lambda)$ in the Bernstein       formula  does not depend on $n$ and is known as the Szeg\"o (or the spectral shift) function. This result can be considered as a far reaching generalization of the asymptotic formula   for the Jacobi polynomials. On the contrary, 
  for singular weights  $\tau(\lambda)$ formula \e{eq:UD1}  is violated but only in neighborhoods of singular points (see, e.g., \cite{M-F}).

  A study of asymptotic behavior  of  orthogonal polynomials defined by their recurrence coefficients $a_{n}$, $b_{n}$ was initiated by P.~Nevai in his book \cite{Nev}.  Under the assumptions  $\{a_{n}-1/2\}\in \ell^1 ({\Bbb Z}_{+})$, $\{b_{n}\}\in \ell^1 ({\Bbb Z}_{+})$ he showed that asymptotics \e{eq:UD1} is true on $\Lambda=(-1,1)$ with $r=0$, $s=1$. 
   The corresponding Jacobi operator is a perturbation
(short-range,   in quantum mechanical terminology)  of the   Jacobi operator  $J^{(0)} :=J_{1/2, 1/2}$ corresponding to the parameters $\alpha=  \beta=1/2$ in \e{eq:Jac1}. It  has the coefficients $a_{n}=1/2$, $b_{n}=0$ for all $n\in {\Bbb Z}_{+}$ and is known as the free discrete Schr\"odinger operator. More general operators $J$ such that $a_{n}\to1/2$, $b_{n}\to 0$ as $n\to\infty$ and
  \begin{equation}
  \sum_{n=0}^\infty (|a_{n+1}-a_{n}| + |b_{n+1}-b_{n}| )<\infty
   \label{eq:LR}\end{equation}
  (long-range perturbations of $J^{(0)}$)
  were considered by A.~Mat\'e, P.~Nevai and V.~Totik in  \cite{Mate}  and in a recent paper \cite{JLR}.  According to the results of these papers an asymptotic behavior of the corresponding orthogonal polynomials $P_{n} (\lambda)$ is always given by formula \e{eq:UD1}  where $r=0$ and $s=1$  but $\Phi_{n}(\lambda)$ is a non-trivial function. 
 Expressions\footnote{Formula (8)  in  \cite{Mate}  is consistent with these expressions but is less explicit.}   (1.13) and (1.22) of
   \cite{JLR}  for the weight $\tau (\lambda)$ and the amplitude factor $\kappa(\lambda)$   confirm identity \e{eq:HXG1}.
   
    \subsection{Freud's weights} 
  
  Starting with the paper \cite{Freud},  exponential  weights
   \begin{equation}
   \tau (\lambda)= k_\beta e^{- |\lambda|^{\beta}}
   \label{eq:Freud}\end{equation}
    where $\beta>0$   and $k_\beta$ is a normalization constant were extensively studied. Note that  the value $\beta=2$  yields the Hermite polynomials.  Given a measure  \e{eq:Lag2}, one constructs its moments, orthonormal polynomials $P_{n} (\lambda)$ and Jacobi coefficients $a_{n}$, $b_{n}$.
      If $ \tau (\lambda)$ is given by \e{eq:Freud} where $\beta\geq  1$, then the minimal Jacobi operator $J_{0}$ with the coefficients  $a_{n}$, $b_{n}$ is essentially adjoint and the spectral measure  of its closure is determined by \e{eq:Freud}.  This is not true for $\beta<1$.
   
        It was shown in \cite{LMS,Magnus} that for weights \e{eq:Freud} and all $\beta>0$,  the off-diagonal recurrence coefficients $a_{n}$ have asymptotics
     \begin{equation}
     a_{n}= \alpha  n^{\ell} (1+ o(1)), \q \ell=1/\beta, \q n\to\infty,
       \label{eq:FreudM}\end{equation}
       with some explicit constant $\alpha=\alpha_{\beta}$; 
of course the diagonal elements $b_{n}=0$.  The asymptotics of orthogonal polynomials $P_{n}(\lambda)$ is given for all $\lambda\in {\Bbb R}$ by  the Plancherel-Rotach formula; see \cite{Nev1}  for $ \beta=4$,  \cite{Rakh}  for   all $ \beta>1$ and \cite{Kriech} for all  $ \beta>0$. Formula \e{eq:PlRot} below states this asymptotics for a fixed $\lambda$  in terms of the coefficients $a_{n}$. 

 
 If   condition \e{eq:FreudM} is satisfied with $ \ell \in (0,1]$ and $b_{n}=0$, then the Carleman condition \e{eq:Carl}  holds and the minimal Jacobi operator $J_{0}$ is essentially self-adjoint.  The  spectrum of its closure $J$   covers the whole axis $\Bbb R$ and it is absolutely continuous so that relation \e{eq:Lag2} holds true with a smooth function $\tau(\lambda)$  for   $\lambda\in \Bbb R$. According to the results of \cite{Jan-Nab} (see also \cite{nCarl}, Sect.~5.2) 
the  corresponding orthonormal polynomials $P_{n} (\lambda)$ have asymptotic behavior
       \begin{equation}
  P_{n} (\lambda)= \pi^{-1/2}  \tau(\lambda)^{-1/2}a_{n}^{-1/2}\cos(\lambda \psi_{n} - n \pi /2+ \d(\lambda))+ o(a_{n}^{-1/2}) 
   \label{eq:PlRot}\end{equation}
  where
   $
 \psi_{n}= 2^{-1}  \sum_{m=0}^{n-1} a_{m}^{-1} 
 $
 and  the phase shift
 $ \d(\lambda) $ does not  depend on $n$.

   It follows that in the case $\ell <1$ formula \e{eq:UD1}  holds true with 
  $r=\ell/2$, $s=1-\ell$, 
  $ \kappa(\lambda)= \pi^{-1/2}    \tau(\lambda)^{-1/2}$,  
  $\omega(\lambda)= 2^{-1} \alpha  (1-\ell)^{-1}\lambda$
  so that  both relations \e{eq:HXG} and \e{eq:HXG1}  are satisfied.   In the case $\ell=1$ we have $\psi_{n}= 2^{-1}\alpha \ln n$  which corresponds to the limit case of \e{eq:HXG} and \e{eq:HXG1} for $s=0$ (see Example~\ref{AOP2x}).

      If condition \e{eq:FreudM}  on $a_{n}$  is satisfied with $\ell > 1$  and all $b_{n}=0$, then, as shown in \cite{nCarl}, the difference equation  \e{eq:Py}  for all $z\in{\Bbb C} $ has solutions $f_{n}^{(+)} (z)$ and $f_{n}^{(-)} (z)
$ (known as   Jost solutions) with asymptotics
          \begin{equation}
f_{n}^{(\pm)} (z)= a_{n}^{-1/2}e^{\pm \pi i n/2}  (1+ o(1)), \q n\to \infty.
\label{eq:JD}\end{equation}
We emphasize  that $z$ is here an arbitrary  complex number and the right-hand side of \e{eq:JD} depends on $z$ only through the remainder $o(1)$. These solutions are linearly independent and belong to $\ell^2 ({\Bbb Z}_{+})$. Therefore the minimal operator $J_{0}$ is not essentially self-adjoint, so that  as shown by R.~ Nevanlinna  (see, e.g., Sect.~7.2 of the book \cite{Schm})   all its  self-adjoint extensions have discrete spectra. Corresponding orthogonal polynomials $P_{n} (z)$ are linear combinations of $f_{n}^{(+)} (z)$ and $f_{n}^{(-)} (z)$ so that an  asymptotic behavior of $P_{n} (z)$ is drastically different from \e{eq:PlRot}.  A value $\ell >1$ corresponds to the case  $\beta<1$ in \e{eq:Freud} which was studied in \cite{Kriech}.  A careful analysis shows that formulas of \cite{Kriech} are consistent with \e{eq:JD}.
   
      \subsection{Jost solutions}

 In the previous subsection,  we discussed a problem with discrete spectrum where  the corresponding difference equation \e{eq:Py} had two solutions with asymptotics  \e{eq:JD}. 
 Actually, it  is rather typical that equation  \e{eq:Py} has two linearly independent solutions   $f_{n}^{(\pm)}  $   with oscillating asymptotics.  Below we  fix some   number $z\in{\Bbb C} $ (not necessarily real).

 \begin{proposition}\label{Jost}
 Suppose that equation  \e{eq:Py}  has two linearly independent solutions  $ f_{n}^{(\pm)}    $
 with asymptotics
  \begin{equation}
f_{n}^{(\pm)}  = v_{n} e^{\pm i\Omega_{n}} (1+ \varepsilon_{n}^{(\pm)}), \q v_{n}>0, \q \Omega_{n}= \ov{\Omega_{n}}.
\label{eq:Jo}\end{equation}

 $1^0$
If $\Omega_{n+1}-\Omega_{n}\to \varpi$ where $\varpi\notin \pi {\Bbb Z}$ and $\varepsilon_{n}^{(\pm)}\to 0$ as $n\to\infty$, then there exists
  \begin{equation}
\lim_{n\to\infty}\big( a_{n} v_{n} v_{n+1}    \big) \neq 0.
\label{eq:Jo2}\end{equation}

 $2^0$
If $\Omega_{n+1}-\Omega_{n}=\pi N + \omega_{n}$ where $N  \in {\Bbb Z} $, $\omega_{n}\to 0$   and $\varepsilon_{n}^{(\pm)}= o (\omega_{n} )$ as $n\to\infty$, then there exists
  \begin{equation}
\lim_{n\to\infty}\big( a_{n} v_{n} v_{n+1}   \omega_{n} \big) \neq 0.
\label{eq:Jo1}\end{equation}
\end{proposition}

      \begin{pf}
       The  Wronskian of the solutions $ \{f_{n}^{(+)} \}$ and $ \{f_{n}^{(-)} \}$  equals
\[
 a_{n}\big( f_{n}^{(+)}  f_{n+1}^{(-)} - f_{n+1}^{(+)}  f_{n}^{(-)} )=2i a_{n} v_{n} v_{n+1}\big(\sin (\Omega_{n}- \Omega_{n+1}) + o(\varepsilon_{n})\big) 
\]
where $\varepsilon_{n}= o(1)$ in case $1^0$ and  $\varepsilon_{n}= o(\omega_{n})$ in case $2^0$.
Since this expression does not depend on $n$ and is not zero, we only have to pass to the limit $n\to\infty$.
   \end{pf}
   
Solutions   $f_{n}^{(\pm)}  $ of   the   Jacobi equation \e{eq:Py} with asymptotics \e{eq:Jo} are known as the Jost solutions. Orthonormal polynomials $P_{n} (\lambda)$ are  their linear combinations   so that
  \begin{equation}
P_{n}  (\lambda) = 2 \kappa(\lambda) v_{n} \cos \big(\Omega_{n}(\lambda) +\d(\lambda)\big)(1+ \varepsilon_{n} (\lambda) )
\label{eq:Jost}\end{equation}
with some coefficients $\kappa(\lambda)$ and $\d(\lambda)$ not depending on $n$.  This formula is consistent with \e{eq:UD1}.

   Let us give some examples where the conditions of Proposition~\ref{Jost} are satisfied.  First we illustrate case  $1^0$. We start with Jacobi operators $J$  whose coefficients $a_{n}\to1/2$, $b_{n}\to 0$ as $n\to\infty$ and obey condition  \e{eq:LR}.  The corresponding orthogonal polynomials contain the classical Jacobi polynomials (see Sect.~2.3). According to Theorem~3.5 in \cite{JLR} we now have
  \[
  \Omega_{n+1} (\lambda) -  \Omega_{n} (\lambda) =\arccos\frac{\lambda-b_{n}}{2a_{n}}, \q \lambda\in (-1,1),
  \]
  so that  $  \varpi=\arccos  \lambda \in (0,\pi)$.  Therefore relation \e{eq:Jo2}  implies that the sequence $v_{n}$ in \e{eq:Jost} has a finite non-zero limit as $n\to \infty$. 
  
   
   Next, we consider recurrence coefficients $a_{n}  $ satisfying   condition \e{eq:FreudM} and $b_{n} = 0$. 
    Now asymptotic formula \e{eq:Jo}  with $  \varpi=-\pi/2$ and $v_{n}= a_{n}^{-1/2} (1+ o(1))$ is true  for $\ell \in (0,1]$ (in particular, the case $\ell  =1/2$ corresponds to  the Hermite polynomials) and all $\lambda\in {\Bbb R}$ according to \e{eq:PlRot}. It is true for $\ell > 1$ and all $z\in {\Bbb C}$ according to \e{eq:JD}. Asymptotic formula for $v_{n}$ is of course consistent with relation \e{eq:Jo2}.

   
   For Laguerre coefficients \e{eq:Lag}, the condition of case $2^0$ is satisfied with $N=1$ and $\omega(\lambda)= 2 \sqrt{\lambda} (\sqrt{n } -\sqrt{n+1 })$ where $\lambda>0$. Therefore relation \e{eq:Jo1}  yields $v_{n}= \nu n^{-1/4}$ for some $\nu >0$ which coincides  with the amplitude factor  in asymptotic formula for the Laguerre polynomials.

\section{Boundedness of integral operators} 

  \subsection{Continuous case}
  
  To motivate our results  on semidiscrete operators  \e{eq:FF2}, let us previously consider the continuous case where the role of $V$ (or of \e{eq:FourF} for $d=1$) is played by the operator
      \begin{equation}
({\sf  V} f) (x)= v (x) \int_{-\infty}^\infty   e^{i\theta(x)\omega (\lambda)}w (\lambda)f(\lambda) d\lambda 
\label{eq:FF1}\end{equation}
with  sufficiently arbitrary functions $\theta(x)$ and $\omega(\lambda)$.  We  suppose that these functions are monotone and tend to $\mp\infty$ (or $\pm\infty$) as $x\to \mp\infty$ and $\lambda\to \mp\infty$. Operators   \e{eq:FF1} can be reduced   to the sandwiched Fourier transforms by   unitary transformations  corresponding to the changes of variables 
$ y= \theta(x), \q \mu=\omega (\lambda)$. 
Indeed, let us set
 \begin{equation}
(F_{1} u)(y)= \theta' (x)^{-1/2} u(x), \q (F_{2} f) (\mu)= \omega' (\lambda)^{-1/2}  f(\lambda).
\label{eq:C4}\end{equation}
Then 
$ F_{1}{\sf V}=  \wt{\sf V} F_{2} $
where the operator $\wt{\sf V}$ acts as
 \[
(\wt{\sf V} \ti{f}) (y)= \ti{v} (y) \int_{-\infty}^\infty e^{i y\mu}  \ti{w}(\mu) \ti{f}(\mu) d\mu
\]
with $\ti{f}= F_{2} f$, 
\[
\ti{v} (y) = \theta' (x)^{-1/2}v(x), \q \ti{w}(\mu) =\omega' (\lambda)^{-1/2} w (\lambda).
\]


Since the operator $\wt{\sf V}$ is bounded if the functions $\ti{v} $ and $\ti{w} $ are bounded, we can state the following elementary result.

 \begin{lemma}\label{CC}
The operator  ${\sf V}: L^2 ({\Bbb R} ) \to L^2 ({\Bbb R})$    is bounded   if  
 \[
(\theta' )^{-1}| v |^2 \in L^\infty ({\Bbb R})\q \mbox{and} \q
(\omega' )^{-1 } |w|^2 \in L^\infty ({\Bbb R}).
\]
 \end{lemma}
 
 Of course the space $L^2 ({\Bbb R} )$ can be replaced here by $ L^2 ({\Bbb R}_{+} )$. For example, the operator  ${\sf V }: L^2 ({\Bbb R}_{+} ) \to L^2 ({\Bbb R}_{+})$   is bounded   if  
 \[
 \theta(x)=x^s, \q |v(x)| \leq C  x^{(s-1)/2}\q\mbox{and}\q \omega(\lambda)=\lambda^t, \q |w(\lambda)|\leq C \lambda^{(t-1)/2}.
 \]
 Note that if $s>1$ or $t>1$, then the functions    $v(x)$  or $w(\lambda)$  may be unbounded at infinity.

The case of operators acting from $ L^2 ({\Bbb R})$ into $ \ell^2 ({\Bbb Z}_{+})$ is essentially more difficult because there is no unitary operator in the space $ \ell^2 ({\Bbb Z}_{+})$ playing the role of the operator $F_{1}$ defined by formula  \e{eq:C4}.  Nevertheless, to a some extent, Lemma~\ref{CC}  remains true in this case -- see Theorem~\ref{BBx}, below.



\subsection{Discrete Fourier transforms}

It is convenient to study boundedness of operators \e{eq:FF2}  in a more general setting. Thus, we consider operators  \e{eq:FourF}  where $\Phi$ is the Fourier transform \e{eq:Four} and   $d M (x)$ 
  is a locally finite measure on ${\Bbb R}^d$.  Let ${\Bbb B}_{r} (x)$ be the ball of radius $r$ centered at a point $x\in {\Bbb R}^d$. We consider measures $dM (x)$  satisfying the condition 
 \begin{equation}
\sup_{x\in {\Bbb R}^d} M({\Bbb B}_{r} (x))<\infty
\label{eq:Four1}\end{equation} 
for some and hence for all $r>0$.

 \begin{theorem}\label{Four}
 If condition \e{eq:Four1} is satisfied and 
  \begin{equation}
|w(\xi)| \leq C (1+ |\xi |^2)^{-\ell/2 }
\label{eq:ww}\end{equation} 
  for some integer $\ell > d/2$, then   operator \e{eq:FourF} is bounded and 
  \begin{equation}
\| A\| \leq C \sup_{x\in {\Bbb R}^d}M ({\Bbb B}_1  (x)).
\label{eq:ww1}\end{equation} 
Conversely, if $w(\xi)$ is the characteristic function of some ball and  the operator $A$  is bounded, then condition \e{eq:Four1} is satisfied.   
\end{theorem}

 \begin{pf}
  Let us consider the lattice of unit cubes $\Pi_{n}$  in ${\Bbb R}^d$. By the Sobolev inequality,   it follows from  condition \e{eq:Four1}  that,  for every $n$,
  \[
  \int_{\Pi_{n}}| u(x)|^2 d M (x)\leq M (\Pi_{n}) \max_{x\in\Pi_{n}} | u(x)|^2  
  \leq C    \int_{\Pi_{n}}\big(| \nabla ^{(\ell)}u(x)|^2 +  | u(x)|^2 \big)d x . 
  \]
  Summing these estimates over all cubes $\Pi_{n}$, we see that
    \begin{equation}
  \|   u\|_{L^2 ({\Bbb R}^d; d M)}^2 \leq C \int_{{\Bbb R}^d}\big(| (\nabla ^{(\ell)}u)(x)|^2 +  | u(x)|^2 \big)d x.
 \label{eq:wsob}\end{equation} 
  It means that  operator  \e{eq:FourF} is  bounded and estimate \e{eq:ww1}  holds.
  
   Conversely, take any not identically zero function $\varphi\in C_{0}^\infty ({\Bbb R}^d)$ such that $w\varphi=\varphi$. There exist a ball ${\Bbb B}_{r_{0}} (x_{0})$  and a number $c>0$ such that the Fourier transform of $\varphi$  satisfies
 \begin{equation}
 |\hat{\varphi}(x)|\geq c>0\q \mbox{for all}\q x\in {\Bbb B}_{r_{0}} (x_{0}).
\label{eq:Four2}\end{equation} 
If condition \e{eq:Four1} is violated, then there exists a sequence $x_{n} $, $x_{n} \to\infty$, such that
 \begin{equation}
  \lim_{n\to\infty } M ({\Bbb B}_{r_{0}} (x_{n})) =\infty. 
\label{eq:Four2a}\end{equation} 
   Set 
   \begin{equation}
 f_{n}  (\xi)=  e^{i(x_{0}-x_{n})\xi}\varphi( \xi).
\label{eq:FGf}\end{equation} 
  Since
  $
  \hat{f}_{n} (x) =   \hat{\varphi} (x-x_{n}  + x_{0}) ,
$
  it follows from \e{eq:Four2}, \e{eq:Four2a}  that
  \begin{equation}
  \| A f_{n}\|_{L^2 ({\Bbb R}^d; d M )}^2 \geq \int_{{\Bbb B}_{r_{0}} (x_{n})}|   \hat{f}_{n} (x) |^2 d M (x)\geq c^2 M({\Bbb B}_{r_{0}} (x_{n}))\to\infty
\label{eq:Four2c}\end{equation}
  as $n\to\infty$. Thus  the operator $A$ is unbounded.
    \end{pf}
    
 Although Theorem~\ref{Four} is   quite elementary, we were unable to localize it in the literature. On the other hand, in the case of absolutely continuous measures, it is equivalent to Theorem~3 in \cite{Y/Osc}. In the general case, the proof is practically the same.
 
 \begin{remark}\label{WSob}
Note that  boundedness of operator \e{eq:FourF} for the characteristic function $w$ of some ball  implies  the same fact   for all functions $w$ satisfying estimate \e{eq:ww}. 
\end{remark}

     \begin{remark}\label{HSH}
Suppose  now that a measure $d M (x)$ is supported by a discrete set $\{x_{1},\ldots, x_{n}, \ldots\}\subset {\Bbb R}^d$  where $n\in{\Bbb Z}_{+}$ or  $n\in{\Bbb Z}$ and set $v_{n}^2=M(\{x_{n}\})$.  Then  condition \e{eq:Four1} means that
 \[
\sup_{x \in {\Bbb R}^d}\sum_{x_{n}\in {\Bbb B}_{r} (x)} v_{n}^2 <\infty.
\]
Note that if the points $x_{n}$ get closer to each other, then   boundedness of  operator \e{eq:FourF}  requires stronger assumptions on decay of $v_{n}$.  Boundedness of $v_{n}$ is always a necessary condition  for boundedness of  operator \e{eq:FourF}.
\end{remark}

    \begin{remark}\label{HSHw}
Condition \e{eq:ww}  in Theorem~\ref{Four} is rather close to optimal.  
Indeed,
even an estimate 
\[
| (A f) (x_n )| \leq C \| f\|, \q \forall f\in C_{0}^\infty ({\Bbb R}), 
\]
for a fixed $x_{n}\in {\Bbb R}^d$
 implies that $ w\in L^2 ({\Bbb R}^d)$.

\end{remark}


  \subsection{One-dimensional case}
  
  Theorem~\ref{Four} applies of course to operators $V: L^2 ({\Bbb R} ) \to \ell^2 ({\Bbb Z}_{+})$ defined by    formula \e{eq:FF2}. This yields   Theorem~\ref{HL}.  Sometimes it may be convenient to state its result in a slightly different form.

 \begin{theorem}\label{FouD}
 Let a sequence $x_{n}\in {\Bbb R}$, $n\in {\Bbb Z}_{+}$,  be such that $x_{n} \to  \infty$ as $n\to   \infty$ and  $x_{n}< x_{n+1}< x_{n}+ \d$ for some $\d>0$ and all $n\in {\Bbb Z}$. 
 If the sequence 
   \begin{equation} 
\sigma_{n}: =\frac{v_{n}^2} {x_{n+1}-x_{n}}
\label{eq:Ms}\end{equation}
 is bounded and condition \e{eq:FF3} on $w(\lambda)$  is satisfied, then the  operator $V$ is bounded. 
Conversely, if $\sigma_{n}\to\infty $ as $n\to\infty$ and  $w(\xi)$ is the characteristic function of some interval, then   $V$ is unbounded. 
\end{theorem}

\begin{example}\label{B1A}
  If   condition \e{eq:rs}  is satisfied,  then $x_{n+1}-x_{n}=s n^{s-1} (1+ o(1))$
    so that $\{\sigma_{n} \}\in \ell ^\infty ({\Bbb Z}_{+})$ and hence the operator $V$ is bounded. This fact can be stated as an inequality  
\begin{equation}
 \sum_{n=1}^\infty n^{s-1} |u(n^s)|^2
 \leq C\int_{0}^\infty (  | u'(x)|^2 +   | u(x)|^2 ) dx.
\label{eq:Sob1A}\end{equation} 
 \end{example}

Let us  also consider the limit case $s=0$, $r=1/2$.
 
  \begin{example}\label{B1a}
  If   
 \[
  x_{n} = \ln n, \q   n >1, \q \mbox{and} \q |v_{n}|
  \leq (1+ | n|)^{-1/2}  ,
\]
then  $x_{n+1}-x_{n}= n^{-1} (1+ o(1))$ so that again $\{\sigma_{n} \}\in \ell ^\infty ({\Bbb Z}_{+})$ and hence the operator $V$ is bounded. This fact can be stated as an
inequality  
\begin{equation}
 \sum_{n=1}^\infty n^{-1} |u(\ln n)|^2
 \leq C\int_{0}^\infty (  | u'(x)|^2 +   | u(x)|^2 ) dx .
\label{eq:Sob1B}\end{equation} 
 \end{example}
 
 Theorem~\ref{FouD}  shows also that inequalities \e{eq:Sob1A}  and \e{eq:Sob1B}  are quite precise, that is, the factors $n^{s-1}$ and $n^{-1}$  cannot  be replaced by $\gamma_{n}n^{s-1}$ and $\gamma_{n}n^{-1}$
 where $\gamma_{n}\to\infty$ as $n\to\infty$.

 \subsection{Additional results}

The first assertion is obvious.
    
     \begin{lemma}\label{HS1}
Let an operator $A$ be defined by relation  \e{eq:FourF}. Then $A$ is Hilbert-Schmidt if and only if $w\in L^2 ({\Bbb R}^d ) $ and $ M ({\Bbb R}^d )<\infty$.
\end{lemma}

Theorem~\ref{Four}  can be supplemented by the following result.

\begin{theorem}\label{FourC}
If condition \e{eq:ww}  is satisfied with an integer $\ell > d/2$ and
 \begin{equation}
\lim_{|x| \to \infty} M  ({\Bbb B}_1 (x))=0,
\label{eq:Four1b}\end{equation}  
  then   the operator $A$  is compact. 
Conversely, if $w(\xi)$ is the characteristic function of some ball and  the operator $A$ is compact, then condition \e{eq:Four1b} is satisfied. 
\end{theorem}

\begin{pf}
Let $\chi_{N} (x)$ be the characteristic function of the ball ${\Bbb B}_{N }(0)$. The operator   $\chi_{N} A$ is Hilbert-Schmidt by Lemma~\ref{HS1}. In view of Theorem~\ref{Four}  it follows from condition \e{eq:Four1b} that
$\|(1- \chi_{N} ) A\| \to 0$ as $N\to\infty$.

Conversely, if condition \e{eq:Four1b} is violated, then 
 there exists a sequence $x_{n} $  such that $|x_{n}| \to\infty$ and
$  M ({\Bbb B}_{r_{0}} (x_{n})) \geq c>0$. Let $f_{n}  (\xi)$ be the same functions \e{eq:FGf} as in Theorem~\ref{Four}.   Similarly to \e{eq:Four2c}, we find  that 
$  \| A f_{n}\|_{L^2 ({\Bbb R}^d; d M)} \geq c>0$.  Since $f_{n}  $ converges weakly to zero,   operator $A$ cannot be compact.  
\end{pf}

Of course an interpolation between Theorem~\ref{Four} and Lemma~\ref{HS1} yields conditions guaranteeing that operators  \e{eq:FourF}  belong to   Schatten-von Neumann classes intermediary between classes of bounded and Hilbert-Schmidt operators, but we will not discuss it here.  Note that very general estimates on singular values of integral operators were obtained 
  by M. Sh. Birman  and M. Z. Solomyak  in the paper \cite{BS2}.


Similarly to Lemma~\ref{CC}, conditions in the variable $\lambda$ in Theorem~\ref{HL}  can be stated in a more general form.
Let us   consider   an  operator ${\bf V}: L^2 ({\Bbb R} ) \to \ell^2 ({\Bbb Z}_{+})$ defined by the  formula 
 \begin{equation}
({\bf V}f)_{n}= v_{n} \int_{-\infty}^\infty   e^{i x_{n} \omega( \lambda )}w (\lambda)f(\lambda) d\lambda
\label{eq:FF2x}\end{equation}
 where  $\omega (\lambda)$ is a  sufficiently  arbitrary function. We suppose that $\omega(\lambda)$ is monotone and $\omega(\lambda)\to\pm \infty$ as $ \lambda\to\pm \infty$. Making
  a change of variables $\mu= \omega (\lambda)$, we see that
  \begin{equation}
({\bf V}f)_{n}=  v_{n} \int_{-\infty}^\infty e^{i x_{n} \mu} \ti{w}(\mu)   \ti{f}(\mu) d\mu
\label{eq:Bb2}\end{equation}
where
$
\ti{f}(\mu) = \omega'(\lambda)^{-1/2}    f(\lambda), \q 
\ti{w}(\mu) = \omega'(\lambda)^{-1/2}    w(\lambda)$.
It follows from Theorem~\ref{HL}  that  operator \e{eq:Bb2}   is bounded if 
\[
| \ti{w}(\mu) |\leq C (1+|\mu|)^{-1}.
\]
Thus Theorem~\ref{HL}  can be supplemented by the following statement.

\begin{theorem}\label{BBx}
Let an operator ${\bf V} : L^2 ({\Bbb R} ) \to \ell^2 ({\Bbb Z}_{+})$ be defined by   formula 
\e{eq:FF2x}  where $\omega' (\lambda)> 0$ and $\omega(\lambda)\to\pm \infty$ as $ \lambda\to\pm \infty$.  Suppose that
  \[
  |w(\lambda)| \leq C \omega' (\lambda)^{1/2}  (1+  |\omega (\lambda)|)^{-1} .
 \]
 Then under   assumption \e{eq:rs} on $x_{n}$ and $v_{n}$  the operator ${\bf V} $ is bounded.
   \end{theorem}
   
   We have chosen the spaces  $ L^2 ({\Bbb R} ) $ and $ \ell^2 ({\Bbb Z}_{+})$ only for definiteness. They can be obviously replaced by $ L^2 ({\Bbb R}_{+} ) $ and $ \ell^2 ({\Bbb Z} )$, respectively.


 \section{Universal relations}
 
  Relation \e{eq:HXG} (but not \e{eq:HXG1}) for the coefficients in asymptotic formula \e{eq:UD1}    can, in principle,  be deduced from the results of Sect.~3. 
  In this section,  we prove both \e{eq:HXG} and  \e{eq:HXG1} by another method, investigating the  time-dependent evolution $\exp(-i\Theta(J) t) f$  as $t\to\infty$ for appropriate functions $\Theta$ of a Jacobi operator $J$.

 \subsection{Stationary phase method} 
 
 Here we consider an integral
  \begin{equation} 
G_{n}  (t)=\int_{\Delta} e^{i \mu x_{n}+i \varphi_{n} (\mu ) -i  \theta (\mu) t}   F(\mu) d\mu 
\label{eq:u}\end{equation}
where $\Delta=(\mu_{1}, \mu_{2})\subset{\Bbb R}$   is a bounded interval, the functions $F\in C_{0}^\infty (\Delta) $, 
 $ \varphi_{n}, \theta\in C^\infty (\Delta) $  and  $\theta' (\mu)> 0$, $\theta'' (\mu) \neq  0$  for $\mu\in \Delta$.  We suppose that $x_{n} <x_{n+1}$, $x_{n}  \to\infty$ and $t\to+\infty$.  With respect to the phases $\varphi_{n} (\mu )$,  we assume\footnote{All estimates of such type are supposed to be uniform in $\mu$ on compact subintervals of $\Delta$.}    that
     \begin{equation} 
  \varphi^{(k)}_{n} (\mu )=o(x_{n}),\q n\to\infty, \q \mu\in \Delta, \q k=1,2,
\label{eq:st}\end{equation}
and $ \varphi^{(k)}_{n} (\mu )=O(x_{n})$  for $k=3,4, \ldots$.
 
Asymptotic behavior of  integral \e{eq:u} for large $t$ and $n$ is determined by stationary  points $\mu$ where
  \begin{equation} 
x_{n}+  \varphi_{n}' (\mu ) -  \theta' (\mu) t =0 .
\label{eq:ust}\end{equation}
A direct integration by parts shows that    integral \e{eq:u}  rapidly tends to zero if the interval $\Delta$ does not contain stationary points.
  
  \begin{lemma}\label{stph1}
     If $x_{n}/ t\not\in \theta' (\Delta)$, then
\begin{equation} 
|G_{n}  (t)|\leq C_k (x_{n}+ t)^{-k} , \q \forall k>0.
\label{eq:u2}\end{equation}
\end{lemma}

Let us introduce the inverse function $h: \theta' (\Delta)\to \Delta$ to $\theta'$ so that    
  \begin{equation} 
\theta' (h(\xi)) =\xi ,\q \xi\in \theta' (\Delta).
\label{eq:ust1}\end{equation}
A next result can easily be obtained by the method of successive approximations.

\begin{lemma}\label{stph}
Let condition \e{eq:st} be satisfied.
 Then for all $\xi_{n}:=x_{n}/t\in \theta' (\Delta)$ and sufficiently large $t$,  equation   \e{eq:ust}
has a unique solution $  \mu_{n}  (t)$  and
 \begin{equation} 
 \mu_{n}( t)= h(\xi_{n})+ o(1), \q t\to+\infty .
\label{eq:st3}\end{equation}
  \end{lemma}

We need the following version of the   stationary phase method.

 \begin{lemma}\label{stph2}
  Let  $\xi_{n} :=x_{n}/ t   \in \theta' (\Delta)$.
Then 
  \begin{equation} 
G_{n}  (t)= \sqrt{2\pi} e^{-\pi i\tau /4} e^{i \Psi_{n} (t) }    t^{-1/2} {\cal F}( x_{n}/t )+  \varrho_{n} (t)\q \mbox{as}  \q t\to+\infty
\label{eq:stph}\end{equation}
where  $\tau= \sgn \theta'' (\mu)$, 
 \begin{equation} 
  {\cal F} ( \xi)= \sqrt{ | h'(\xi)| }  F( h (\xi)),
\label{eq:stFG}\end{equation}
 \begin{equation} 
\Psi_{n}  (t)=   \mu_{n} (t) x_{n}+  \varphi_{n} ( \mu_{n} (t)) -  \theta ( \mu_{n} (t)) t 
\label{eq:stph1}\end{equation}
and
 \begin{equation} 
\sup_{  x_{n}/t \in \theta' (\Delta)  }| \varrho_{n} (t)|= o(t^{-1/2}).
\label{eq:u1}\end{equation}
\end{lemma}

{\it Sketch of the proof.}  Applying the   stationary phase formula to integral \e{eq:u}, we see that
 \begin{equation} 
G_{n}  (t)= \sqrt{2\pi}   e^{-\pi i\tau/4} e^{i \Psi_{n} (t) }   |\theta''  ( \mu_{n} (t)) t - \varphi_{n}''( \mu_{n} (t))|^{-1/2} F ( \mu_{n} (t)) )+  \varrho_{n} (t) 
\label{eq:Stph}\end{equation}
where $r_{n} (t) $ satisfies \e{eq:u1}
as $t\to+\infty$.  According to \e{eq:st3}  we can here replace  $\mu_{n} (t)$ by $ h (\xi_{n})$ and according to \e{eq:st}  we can neglect $\varphi_{n}''(\mu_{n} (t)  )$.
Differentiating \e{eq:ust1}, we see that $\theta'' (h(\xi)) h'(\xi)=1$. Therefore using notation \e{eq:stFG}, we can rewrite \e{eq:Stph} as \e{eq:stph}.

 \subsection{Asymptotic evolution} 
 
  Let $\Lambda$ be an interval of
  the absolutely continuous spectrum of a Jacobi operator $J$  so that its  spectral measure satisfies \e{eq:Lag2}  where we suppose that the weight $\tau\in C^\infty (\Lambda)$ and $\tau(\lambda)>0 $  for $\lambda\in\Lambda$. With respect to   the corresponding orthonormal polynomials $P_{n} (\lambda)$, we accept an assumption somewhat more general than \e{eq:UD1}.
  
  \begin{assumption}\label{ASY}
 As $\lambda\in \Lambda$ and $n\to\infty$, an asymptotic formula  holds
  \begin{equation} 
P_{n}  (\lambda) = 2 \kappa   (\lambda)v_{n} \cos \big(   \omega (\lambda) x_{n} + \Phi_{n}(\lambda)\big)+ \epsilon_{n}  (\lambda)
\label{eq:HX1M}\end{equation}
where  $x_{n+1}> x_{n}$,  $x_{n}\to \infty$   as $n\to\infty$, $v_{n}>0$ and
  \begin{equation} 
\sum_{n=1}^\infty  v_{n}^2 x_{n}^{-2k}< \infty
\label{eq:Nvx}\end{equation}
for some $k>0$. We suppose that
$\omega, \kappa \in C^\infty (\Lambda)$,   $\kappa(\lambda)>0$, $\omega' (\lambda) > 0$,   
  \[
    \Phi^{(k)}_{n} (\lambda )=o(x_{n}),\q n\to\infty, \q \lambda\in \Lambda, \q k=1,2,
\]
 and $  \Phi^{(k)}_{n} (\lambda )=O(x_{n})$ for all $k \geq 3$.
The remainder in \e{eq:HX1M} satisfies a condition $\sum_{n} |\epsilon_{n}  (\lambda)|^2 <\infty$
uniformly on compact subintervals of $\Lambda$.
\end{assumption}


 Let $f \in C_{0}^\infty (\Lambda)$ and $u= Uf \in L^2 (\Lambda)$  so that $u= E (\Lambda) u$.  It follows from formulas \e{eq:UFf} and \e{eq:HX1M}  that, for an arbitrary function $\Theta\in L^\infty (\Lambda)$, 
    \begin{equation} 
    ( e^{-i \Theta (J) t} u)_{n}= 2 v_{n}\int_{\Lambda}  \cos \big(   \omega (\lambda) x_{n} + \Phi_{n}(\lambda)\big) 
 e^{-i \Theta (\lambda) t} \varkappa (\lambda)  f (\lambda) d\lambda  +\varrho_{n} (t)
\label{eq:N}\end{equation}
where    $   \varkappa (\lambda) =\kappa(\lambda)\sqrt{\tau(\lambda)}$ and 
   \begin{equation} 
   \varrho_{n} (t) =   \int_{\Lambda} \epsilon_{n}(\lambda) e^{-i \Theta(\lambda) t}\sqrt{\tau(\lambda)}  f(\lambda) d\lambda.
   \label{eq:fg}\end{equation}
    Let us
   find an asymptotic behavior in $\ell^2 ({\Bbb Z}_{+})$ of sequence \e{eq:N} as $t\to\infty$. 
   The last   term on  the right 
    is negligible. Indeed, the operator $e^{-i \Theta(J) t}$ tends weakly to zero as $t\to\infty$ if $\Theta' (\lambda)\neq 0$ for $\lambda\in \Lambda$. The operator $R: L^2 (\Lambda)\to \ell^2 ({\Bbb Z}_{+})$ defined by 
$
(R f )_{n}= \int_{\Lambda_{c}} \epsilon_{n}(\lambda) \tau(\lambda) f(\lambda) d\lambda
$
is in the Hilbert-Schmidt class for any compact subinterval $\Lambda_{c}$ of $\Lambda$. It follows that 
\begin{equation} 
\lim_{t\to\infty}\sum_{n =0}^\infty  | \varrho_{n} (t)|^2=0.
\label{eq:NNx}\end{equation} 


Assume now that
\begin{equation} 
\Theta (\lambda) =\theta (\omega(\lambda))
\label{eq:NN}\end{equation}
for some function $\theta\in C^\infty (\Delta)$ where   $\theta' (\mu)>0 $ for $\mu\in \Delta=\omega (\Lambda)$.
 Making the change of variables 
\begin{equation} 
\mu= \omega(\lambda),\q \phi_{n}(\mu)=   \Phi_{n}(\lambda) , \q F(\mu)=\omega' (\lambda)^{-1}\varkappa (\lambda)  f(\lambda)
\label{eq:N1}\end{equation}
in the first   term in  the right-hand side of \e{eq:N},  we find that it equals
 \begin{equation} 
 v_{n}     \int_{\Delta}  
 e^{i \mu x_{n} + i \phi_{n}(\mu)-i \theta (\mu) t} F(\mu) d \mu +   v_{n}  \int_{\Delta}  
 e^{-i \mu x_{n} - i \phi_{n}(\mu)-i \theta (\mu) t} F(\mu) d \mu.
\label{eq:N11}\end{equation}
In the region $x_{n} /t  \in \theta' (\Delta)$, 
  asymptotics of the first integral in  \e{eq:N11} is given by  Lemma~\ref{stph2}. 
A  direct integration by parts  shows that for $x_{n} /t \not \in \theta' (\Delta)$, this integral  satisfies estimate \e{eq:u2}  so that under   assumption \e{eq:Nvx},  this region of $n$ can be neglected. Similarly, integrating by parts we see that
the second integral  satisfies estimate \e{eq:u2}  for all $n$.  Therefore    the second term in  \e{eq:N11} tends to zero  in $\ell^2 ({\Bbb Z}_{+})$ as $t\to\infty$.  This leads to  the following result.

\begin{theorem}\label{ev}
 Let Assumptions~\ref{ASY}   be satisfied, and let  $f\in C_{0}^\infty (\Lambda) $, $u=Uf$.   Define the function $\Theta(\lambda)$ by equality \e{eq:NN}, the functions $\varphi_{n}  (\mu)$, $F (\mu)$ -- by formula \e{eq:N1},  ${\cal F}  (\xi) $ -- by \e{eq:stFG} and $\Psi_{n} (t)$ -- by  \e{eq:stph1}. 
    Then, as $t\to+\infty$, 
    \begin{equation} 
    ( e^{-i \Theta (J) t} u)_{n}= \sqrt{2\pi} v_{n} e^{-\pi i\tau/4}  e^{i\Psi_{n}(t)} t^{-1/2} {\cal F}(x_{n}/t)\mathbbm{1} ( x_{n}/ t  \in {\cal I})   + \varrho_{n}  (t), 
\label{eq:NX}\end{equation}
where $ {\cal I}=\theta' (\omega (\Lambda)) $, $\tau= \sgn \theta''$ and the remainder $\varrho_{n}  (t)$ satisfies \e{eq:NNx}.
\end{theorem}

According to \e{eq:NX} the sequence $ (e^{-i \Theta (J) t } u)_{n}$ ``lives" in the region of $n$ such that $x_{n} / t \in {\cal I}$. Below, the phase factor $e^{i\Psi_{n}(t)} $ is not very important. Since 
  \begin{equation} 
\| e^{-i \Theta (J) t } u \| =\| u\| =\| Uf \| = \| f\|,
\label{eq:UX}\end{equation}
we can  state   

\begin{corollary}\label{ev1}
Under the assumptions of Theorem~\ref{ev}  we have
  \begin{equation} 
2\pi \lim_{t\to\infty}  t^{-1}\! \!\!\sum_{n : x_{n} /t  \in {\cal I}}  v_{n}^2 \,  | {\cal F}(x_{n}/t)|^2= \| f\|^2.
\label{eq:NX1}\end{equation}
\end{corollary}

\begin{remark} 
It  suffices to make   assumptions on only one derivative  of functions  $\tau, \kappa$ and two  derivatives of functions   $\omega, \Phi_{n}$. However in this case one has to demand  \e{eq:Nvx} for $k=1$.
\end{remark}

\begin{remark} 
The stationary phase method does not necessarily require that the function $\Theta(\lambda)$ be defined by formula \e{eq:NN}.  If $\omega'' (\lambda)\neq 0$, we can also take $\Theta(\lambda)=\lambda$ and if $\omega(\lambda) = \omega_{0}\lambda$, any function $\Theta(\lambda)$ with $\Theta''(\lambda)\neq 0$  is admissible. Asymptotic formula \e{eq:NX}  remains essentially the same in all  these cases. On the contrary  if 
$\omega(\lambda)=\omega_{0}\lambda$ and $\Theta(\lambda)=\Theta_{0}\lambda$, then the evolution $e^{-i \Theta (J) t} u$ reduces to shifts; see Appendix~ B.
\end{remark}

 \subsection{Main result} 
 
 Relations  \e{eq:HXG} and \e{eq:HXG1}  are deduced from formula \e{eq:NX1}. To that end, we need
 the following auxiliary result.  It is a direct consequence of the definition of an integral of a continuous function as a limit of the corresponding Riemann sums.
 
 \begin{lemma}\label{INT}
 Suppose that ${\cal I}\subset{\Bbb R}$ is a bounded interval and a function ${\cal F}\in C({\cal I})$. Let 
 a sequence $X_{n} (t)\in {\Bbb R}$, $n\in {\Bbb Z}_{+}$,  be such that $X_{n+1} (t)> X_{n }(t)$, $X_{n} (t)\to   \infty$ as $n\to   \infty$ for all $t$ and
    \begin{equation}
\lim_{t\to\infty} \max_{n : X_{n} (t)\in {\cal I}}  (  X_{n+1} (t)- X_{n} (t))=0.
\label{eq:INT1}\end{equation} 
 Then 
    \begin{equation}
 \lim_{t\to\infty }\sum_{n : X_{n}(t) \in {\cal I} } (X_{n+1}  (t)- X_{n}(t)) |{\cal F}(X_{n}(t))|^2=\int_{{\cal I}} |{\cal F}(\xi) |^2  d \xi .
\label{eq:INT}\end{equation} 
\end{lemma}

  Now we can state our main result.
 
  \begin{theorem}\label{ev2}
   Let Assumption~\ref{ASY} be satisfied. Suppose additionally that
  \begin{equation}
\lim_{n\to\infty} (x_{n+1}-x_{n}) /x_{n}=0.
\label{eq:INT2}\end{equation} 
Define a sequence $\sigma_{n}$ by equality \e{eq:Ms}
and assume that a   limit $($possibly,  infinite$)$
 \begin{equation} 
\lim_{n\to\infty}\sigma_{n}=\sigma
\label{eq:M}\end{equation}
exists. Then necessarily $0<\sigma<\infty$ 
and
 \begin{equation} 
2 \pi \sigma\tau(\lambda) \kappa^2 (\lambda) =    \omega' (\lambda) .
\label{eq:M1}\end{equation}
 \end{theorem}

  \begin{pf}
  Let us set ${\cal I}= \theta' (\omega (\Lambda))$, $X_{n}(t)=x_{n}/t$ and observe that the sequence
  \[
  X_{n+1}(t) -X_{n}(t)=  ((x_{n+1}-x_{n}) /x_{n})   X_{n}(t)
  \]
  satisfies condition \e{eq:INT1}.
 Indeed,  the first factor tends to zero in view of \e{eq:INT2} and the second one is bounded if $X_{n}(t)\in {\cal I}$.
 Note also that $n\to\infty$ if $X_{n}(t)\in {\cal I}$ and $t\to\infty$.
 Applying relation \e{eq:INT}  to a function $ {\cal F}(\xi)$ defined by equalities \e{eq:stFG} and \e{eq:N1} and using condition  \e{eq:M}, we find that 
   \begin{multline} 
 \lim_{t\to\infty }\sum_{n : X_{n}(t) \in {\cal I} } \sigma_{n}(X_{n+1}  (t)- X_{n}(t)) |{\cal F}(X_{n}(t))|^2= \sigma \int_{{\cal I}} | {\cal F}(\xi) |^2  d \xi 
\\
 =  \sigma \int_{{\cal I}} | h'(\xi)| | F(h (\xi)) |^2  d \xi 
= \sigma \int_{\omega (\Lambda)}   | F(\mu) |^2  d\mu = \sigma \int_{ \Lambda}  \omega' (\lambda)^{-1}\varkappa (\lambda)^2   | f(\lambda) |^2  d\lambda. 
\label{eq:NXX}  \end{multline} 
Since 
$
t^{-1}v_{n}^2 =\sigma_{n}(X_{n+1}  (t)- X_{n}(t)),
$
it follows from Corollary~\ref{ev1} that the left-hand side of  \e{eq:NXX}  equals $(2\pi)^{-1}\| f\|^2$.
Therefore relation \e{eq:NXX}  ensures that
 \begin{equation} 
 \int_{\Lambda}    |f(\lambda) |^2  d\lambda
= 2\pi \sigma \int_{\Lambda}  \omega'(\lambda)^{-1}\varkappa^2(\lambda) |f(\lambda) |^2  d\lambda 
\label{eq:MK}\end{equation}
for an arbitrary function $f\in C_{0}^\infty (\Lambda)$.  Thus,  $\sigma$ is necessary finite, $\sigma\neq 0$ and 
  $2\pi \sigma  \omega'(\lambda)^{-1}\varkappa^2(\lambda)=1$. This yields \e{eq:M1}.
    \end{pf}
    

     \begin{remark}\label{OPx}
     Condition  \e{eq:INT2} imposes a  mild restriction on the growth of $x_{n}$ as $n\to\infty$. However it is satisfied for $x_{n }= n^s$ with an arbitrary $s>0$.
 \end{remark}
 
 \subsection{Examples and comments}
 
 The first example generalizes the results stated in Sect.~1.1. The second corresponds to the limit case $s=0$.
 

\begin{example}\label{AOP2}
Suppose that asymptotic formula \e{eq:HX1M} holds with
 \begin{equation} 
 x_{n }= n^s (\ln n )^p , \q  v_{n}=n^{-r} ( \ln n)^q , \q r \geq 0,
\label{eq:M2x}\end{equation}
where $s>0$ and  $p,q \in {\Bbb R}$. Then 
  number \e{eq:Ms}  equals
\[
\sigma_{n}  =s^{-1} n^{-2r -s+1} (\ln n)^{2q-p} (1+ o(1)), \q n\to\infty,
\]
so that limit \e{eq:M} exists. It is finite and $\sigma\neq 0$ if the conditions \e{eq:HXG} and $2q=p$   are satisfied. Relation \e{eq:M1} holds true with $\sigma=s^{-1}$. 
\end{example}

   \begin{example}\label{AOP2x}
Suppose that \e{eq:HX1M} and \e{eq:M2x} hold with $s=0$.  Then we have
\[
\sigma_{n}  =p^{-1} n^{-2r +1} (\ln n)^{2q-p+1} (1+ o(1)), \q n\to\infty,
\]
and hence  limit \e{eq:M}  again  exists. It is finite and $\sigma\neq 0$ if the conditions $r=1/2$ and $2q=p-1$   are satisfied. Relation \e{eq:M1} holds true with $\sigma=p^{-1}$. 
 \end{example}

   Finally, we note that both the conditions and the statements of Theorem~\ref{ev2}  are   quite different from those of Proposition~\ref{Jost}. In particular, the latter   imposes assumptions on the phase $\Omega_{n} $  while  Theorem~\ref{ev2} requires  conditions on its derivative   $\Omega_{n}' (\lambda) $. Note  also that
   Theorem~\ref{ev2} is specific for the absolutely continuous spectrum of $J$.
  




 \appendix
\section{An elementary inequality }


As a by-product of our considerations,  we obtain a simple inequality for functions in the Sobolev space ${\sf H}^1({\Bbb R})$.  First, we note that
  \begin{equation}
| u(0)|^2 \leq c_{0}   \int_{0}^{1}( | u'(x)|^2 +  | u(x)|^2 ) d x 
\label{eq:Sob}\end{equation} 
where we can choose $c_{0} =2(\sqrt{5}-1)^{-1}$.  We check the following result.

 \begin{proposition}\label{Sob1}
 Let a sequence $x_{n}\in {\Bbb R}$, $n\in {\Bbb Z}$,  be such that $x_{n} \to \pm \infty$ as $n\to \pm \infty$ and  $x_{n}< x_{n+1}< x_{n}+ \d$ for some $\d>0$ and all $n\in {\Bbb Z}$. Then
   \begin{equation}
 \sum_{n\in {\Bbb Z}} (x_{n+1}-x_{n}) |u(x_{n})|^2
 \leq c_{0}\max\{1,  \d^2 \}\int_{-\infty}^\infty (  | u'(x)|^2 +   | u(x)|^2 ) dx
\label{eq:Sob1}\end{equation} 
for all $u\in{\sf H}^1({\Bbb R})$.
\end{proposition}

 \begin{pf}
 Applying estimate \e{eq:Sob} to the function $u_{\varepsilon}  (x)= u(\varepsilon x)$ and making the change of variables $y=\varepsilon x$, we see that
  \[
| u(0)|^2 \leq c_{0}   \int_{0}^{1}(\varepsilon^2 | u '(\varepsilon x)|^2 +  | u (\varepsilon x)|^2 ) dx 
=c_{0}   \int_{0}^{\varepsilon}(\varepsilon | u'(y)|^2 + \varepsilon^{-1} | u (y)|^2 ) dy 
\]
whence
\begin{equation}
\varepsilon | u(x_{n})|^2 \leq c_{0}   \int_{x_{n}}^{x_{n} +\varepsilon}(\varepsilon^2 | u'(y)|^2 +  | u (y)|^2 ) d y . 
\label{eq:Sob3}\end{equation} 
Let us set here $\varepsilon =x_{n+1} -x_{n}$. Then $\varepsilon\leq\d$  and \e{eq:Sob3} yields
\[
(x_{n+1} -x_{n}) | u(x_{n})|^2 \leq c_{0} \max\{1,  \d^2 \}  \int_{x_{n}}^{x_{n+1}  }(  | u'(y)|^2 +  | u (y)|^2 ) d y . 
\]
It remains to take the sum of these estimates over all $n\in {\Bbb Z}$. 
  \end{pf}

  Inequality \e{eq:Sob1} is   quite elementary but, surprisingly, we were unable to find it in the literature.  Of course
 \e{eq:Sob1} implies again inequalities  \e{eq:Sob1A} and \e{eq:Sob1B}.

  
   \section{Dispersionless evolution}
  
  If $x_{n}= o(n)$ as $n\to\infty$ (which excludes the linear growth of $x_{n} $)   and the phases $\Phi_{n} $ do not depend on $\lambda$, we can dispense with the stationary phase method.     Here we consider  evolution \e{eq:N}  for the case $\Theta(\lambda)=\omega(\lambda)$ excluded in Sect.~4.2.  Making the change of variables \e{eq:N1}, we can rewrite \e{eq:N} as
  \begin{equation} 
    ( e^{-i \omega (J) t} u)_{n}= \sqrt{2\pi} v_{n} \big(e^{i\Phi_{n}} \wh{F} (x_{n}-t) + e^{-i\Phi_{n}} \wh{F} (-x_{n}-t)\big)+ \varrho_{n}  (t)
\label{eq:Nx}\end{equation}
where $\wh{F} (x) $ is the Fourier transform of $F(\mu)$. The remainder $ \varrho_{n} (t)$ is given by \e{eq:fg}; it is negligible according to \e{eq:NNx}. Since $\wh{F} (-x )=O(x^{-k})$ as $x\to+ \infty$ for all $k$, the term with $\wh{F} (- x_{n}-t)$ in \e{eq:Nx} is also negligible. Thus instead of Theorem~\ref{ev}  we now have

\begin{theorem}\label{evD}
 Let Assumption~\ref{ASY}   be satisfied, and let  $f\in C_{0}^\infty (\Lambda) $.  Then 
   \begin{equation} 
    ( e^{-i \omega (J) t} u)_{n}= \sqrt{2\pi} v_{n}  e^{i\Phi_{n}} \wh{F} (x_{n}-t)   + \varrho_{n}  (t)
\label{eq:NxC}\end{equation}
 where the remainder $\varrho_{n}  (t)$ obeys condition  \e{eq:NNx}.   
 \end{theorem}

It follows from \e{eq:NxC} that
 \[
\|  e^{-i\omega (J) t}  u\|^2= 2\pi \sum_{n\in {\Bbb Z}_{+}}v_{n}^2 | \wh{F} (x_{n}-t)|^2+ o(1), \q t\to\infty,
\]
and hence according to \e{eq:UX} (where $\Theta=\omega$)
 \begin{equation} 
2\pi \lim_{t\to\infty } \sum_{n\in {\Bbb Z}_{+}}v_{n}^2 | \wh{F} (x_{n}-t)|^2 = \| f\|^2.
\label{eq:N2U}\end{equation}

     Note that formula \e{eq:NX} is not true now because 
if $\Theta(\lambda)=\omega(\lambda)$, then
$\theta(\lambda)=\lambda$ so that $\theta'' (\lambda)=0$.

On the other hand, one can calculate the left-hand side of \e{eq:N2U} using Lemma~\ref{INT}. Indeed, let us   set
  $X_{n}  (t)=x_{n } -t$ and observe that  $v_{n}^2= \sigma_{n} (X_{n+1} (t) -X_{n} (t))$. If the sequence $\sigma_{n} $ has a finite limit $\sigma$, then according to  Proposition~\ref{Sob1}, for each $R>0$,   the sum in \e{eq:N2U} over $n$ such that $|X_n (t)|\geq  R$ is estimated by 
  $C \int_{| x |\geq R} (|\wh{F}'(x)|^2 + |\wh{F}( x )|^2)d x$. In view of  \e{eq:INT}  the sum in \e{eq:N2U} over $n$ such that $|X_n (t)| <  R$  converges to $\sigma \int_{-R}^R | \wh{F}(x) |^2 d x$ as $t\to\infty$. It follows that the left-hand side of \e{eq:N2U}  equals $   2\pi \sigma \|  \wh{F}\|^2 = 2\pi \sigma \|  F\|^2$ so that $    2\pi \sigma \|  F\|^2 = \| f\|^2   $. This yields again equality \e{eq:MK} whence
   $\sigma>0$ and relation   \e{eq:M1}      follows.


 \end{document}